\documentclass{article}

\usepackage[final]{nips_2018}
\usepackage{graphicx, multicol, latexsym, amsmath, amssymb}
\usepackage[utf8]{inputenc} 
\usepackage[T1]{fontenc}    
\usepackage{hyperref}       
\usepackage{url}            
\usepackage{booktabs}       
\usepackage{xspace} 
\usepackage{amsfonts}       
\usepackage{nicefrac}       
\usepackage{microtype}      
\usepackage{amsmath}
\usepackage{subfigure}
\usepackage{comment}
\newtheorem{theorem}{Theorem}[section]
\newtheorem{definition}{Definition}[section]
\newcommand{\norm}[1]{\left\lVert#1\right\rVert}
\newcommand{\delf}{DELF\xspace}
\DeclareMathOperator*{\argmin}{arg\,min}
\DeclareMathOperator*{\expectation}{\mathbb{E}}
\title{Deep Lyapunov Function: Automatic Stability Analysis for Dynamical Systems}

\author{Arash Mehrjou\\
Department of Empirical Inference\\
Max Planck Institute for Intelligent Systems\\
\texttt{arash.mehrjou@tuebingen.mpg.de}\\
\And
Bernhard Sch\"{o}lkopf\\
Department of Empirical Inference\\
Max Planck Institute for Intelligent Systems\\
\texttt{bs@tuebingen.mpg.de} \\
}

\begin{document}

\maketitle

\begin{abstract}
Stability analysis plays a crucial role in studying the behavior of dynamical systems with theoretical and engineering applications. Among various kinds of stability, the stability of equilibrium points is of the greatest importance which is mainly studied by Lyapunov's stability theory. This theory requires finding a function with specified properties. Except for a few simple examples, there is no straightforward constructive algorithm to find a Lyapunov function for an arbitrary dynamical system. The goal of this work is proposing a simple yet effective way to approximate this function using deep learning tools. 
\end{abstract}

\section{Introduction}
Dynamical systems are often described by a set of coupled differential equations:
\begin{equation}
\label{eq:statespace}
\begin{bmatrix}
       \dot{x}_1\\
       \dot{x}_2\\
       \vdots\\
       \dot{x}_n
\end{bmatrix} = 
\begin{bmatrix}
       f_1(x_1,x_2,\hdots,x_n,u_1,u_2,\hdots,u_m,t)\\
       f_1(x_1,x_2,\hdots,x_n,u_1,u_2,\hdots,u_m,t)\\
       \vdots\\
       f_1(x_1,x_2,\hdots,x_n,u_1,u_2,\hdots,u_m,t)\\
\end{bmatrix}
\end{equation}
which can be more compactly represented by $\dot{x}(t)=f(x,u,t)$ where $x, f \in R^n$ and $u \in R^m$.
The state vector $x(t)$ fully determines the system at each time instant. Input vector $u(t)$ forces the system to achieve a determined goal. In this work, we assume the system is closed-loop, i.e. $u(t)=u(x(t))$. Therefore, we drop $u(t)$ from the equations. In addition, we only consider time-invariant systems, meaning that the explicit dependence of the right-hand side of~\eqref{eq:statespace} on time $t$ is dropped. Systems with these properties are called autonomous systems and simply represented by
\begin{equation}
\label{eq:statespace_compact}
\dot{x}=f(x)
\end{equation}
We also dropped time argument of $x(t)$ for writing convenience. Existence and uniqueness of the solution of this system is guaranteed if $f:D\to R^n$ is a locally Lipschitz map from a domain $D\subset R^n$ into $R^n$. Suppose $\bar{x}\in D$ is an equilibrium point of system \eqref{eq:statespace_compact}; that is, $f(\bar{x})=0$. Stability analysis concerns the behavior of the system in a vicinity of its equilibrium point $\bar{x}$. Without loss of generality, we can always change the coordinates such that the equilibrium point of the system sits in the origin of $R^n$~\citep{khalil1996noninear}. Therefore, from now on we assume $\bar{x}=0$. In the following, we briefly provide the definition of different stability conditions about the equilibrium point and afterwards present one of the most important theorems in stability analysis of dynamical systems.
\begin{definition}
The equilibrium point $x=0$ of \eqref{eq:statespace_compact} is
\begin{itemize}
\item \textbf{stable} if for each $\epsilon>0$, there is $\delta=\delta(\epsilon)>0$ such that
\begin{equation}
\label{eq:def_stable}
\norm{x(0)}<\delta \implies \norm{x(t)}<\epsilon, \quad \forall t \leq 0
\end{equation}
\item \textbf{unstable} if it is not stable.
\item \textbf{asymptotically stable} if it is stable and $\delta$ can be chosen such that 
\begin{equation}
\norm{x(0)}<\delta \implies \lim_{t\to \infty}x(t)=0
\end{equation}
\end{itemize}
\label{def:stablility}
\end{definition}
Given these definitions, in the following, Lyapunov's first stability theorem provides a way to determine the stability condition of the equilibrium point.
\begin{theorem}
Let $x=0$ be an equilibrium point for system \eqref{eq:statespace_compact} and $D\subset R^n$ be a domain containing $x=0$. Let $V:D\to R$ be a \emph{continuously differentiable function} such that 

\begin{align}
V(0)&=0 \quad {\rm and}\quad V(x)>0\quad {\rm in}\quad D-\{0\}\label{eq:lyapanov_psd}\\
\dot{V}(x)&\leq0 \quad {\rm in}\label{eq:lyapanov_nsd}\quad D
\end{align}
Then, $x=0$ is \textbf{stable}. Moreover, if 

\begin{equation}
\dot{V}(x)<0\quad{\rm in}\quad D-\{0\}
\label{eq:lyapanov_nd}
\end{equation}

then $x=0$ is \textbf{asymptotically stable}.
\label{theo:Lyapanov}
\end{theorem}

A continuously differentiable function $V(x)$ that satisfies both \eqref{eq:lyapanov_psd} and \eqref{eq:lyapanov_nsd} is called a \emph{Lyapanov function}. An intuitive description of equations \eqref{eq:lyapanov_psd} and \eqref{eq:lyapanov_nsd} is that when the state trajectory of the system enters the set $\Omega_c=\{x\in R^n|V(x)\leq c\}$, it never comes out of it. When $\dot{V}<0$ as in \eqref{eq:lyapanov_nd}, the set $\Omega_c$ shrinks and the trajectory will approach the origin $x=0$. Function $V(x)$ satisfying \eqref{eq:lyapanov_psd} is said to be \emph{positive definite}. Instead, when a weaker condition $V(x)\geq 0$ is satisfied, it is called \emph{positive semidefinite}. Likewise, $V(x)$ satisfying \eqref{eq:lyapanov_nsd} is said to be \emph{negative definite} and is known as \emph{negative semidefinite} when $V(x)\leq 0$. If none of these conditions hold, the function $V(x)$ is said to be \emph{indefinite}.

In practice, the message of theorem~\ref{theo:Lyapanov} is finding a function $V(x)$ which satisfies the required conditions~\eqref{eq:lyapanov_psd} and one of \eqref{eq:lyapanov_nsd} or \eqref{eq:lyapanov_nd}. There is no general constructive way for finding this function~\citep{lyapunov1992general}. The path is though more clear in some occasions~\citep{hafstein2007algorithm}. For example, when the system of interest is a physical system, a function $E(x)$ that characterize the energy of the system can be a natural candidate for a Lyapunov function. Notice that even in this case, there is no guarantee that the energy function is the best Lyapunov function for characterizing the stability properties of the system. Search for a better Lyapunov function is motivated by its use in downstream tasks including estimating the domain of attraction, controller design strategies such as Lyapunov redesign method~\citep{nevsic2005lyapunov}, and achieving better understanding of the qualitative and quantitative behavior of the dynamical system~\citep{isidori2013nonlinear}. Moreover, apart from engineering applications, Lyapunov analysis is widely used in studying the stability and convergence of iterative optimization and learning algorithms~\citep{wilson2016lyapunov}.

In this work, a simple but effective method for searching for a Lyapunov function is proposed. To this end, we take advantage of a well-known property of multilayer perceptrons(MLP) being \emph{universal function approximator}~\citep{csaji2001approximation} combined with stochastic gradient descent (SGD) as a generic optimization method which is widely used in deep learning~\citep{robbins1985stochastic, goodfellow2016deep}. 

\section{Proposed Method: \delf}
DEep Lyapunov Function (\delf) is proposed as a method to automate finding a Lyapunov function for the dynamical system \eqref{eq:statespace_compact}. We parameterize a  scalar function $\hat{V}(x;\theta): R^n \to R$ by a deep neural network. The idea is to find $\theta$ such that \eqref{eq:lyapanov_psd} and \eqref{eq:lyapanov_nd} are satisfied for $\hat{V}$. We take an empirical approach to check this satisfaction which of course is not a rigorous mathematical guarantee and its success depends on the effectiveness of SGD in finding the parameters of MLPs. Because the desired value of $\theta$ must make $V$ positive definite and $\dot{V}=\nabla_x\hat{V}(x;\theta)^Tf(x)$ negative definite, we propose the following loss function

\begin{equation}
\mathcal{L}(x;\theta)=h_1(\hat{V}(x;\theta))+h_2(\nabla_x\hat{V}(x;\theta)^Tf(x))
\label{eq:loss}
\end{equation}

with the following description of its components: The first term corresponds to positive definiteness of $V$ and the second term corresponds to negative definiteness of $\dot{V}(x)=\nabla_x\hat{V}(x;\theta)^Tf(x)$. Since we wish $V$ to be positive definite, negative values of $V$ must be punished. To avoid collapsing $V$ on zero, a margin $m_1$ is introduced and $h_1(.)$ is defined as

\begin{equation}
h_1(V) =
\left\{
	\begin{array}{ll}
		0  & \mbox{if } V > m_1 \\
		-V+m_1 & \mbox{if } V \leq m_1.
	\end{array}
\right.
\label{eq:h1}
\end{equation}

Similarly, the second term of the loss function pushes $\dot{V}$ towards negative values less than a margin $-m_2$ when function $h_2(y)$ is defined as
\begin{equation}
h_2(\dot{V}) =
\left\{
	\begin{array}{ll}
		\dot{V}+m_2  & \mbox{if } \dot{V} > - m_2 \\
		0 & \mbox{if } \dot{V} \leq - m_2.
	\end{array}
\right.
\label{eq:h2}
\end{equation}

Notice that $m_1, m_2 \geq 0$. In general, there is no need to define $h_1$ and $h_2$ similarily or set the margines $m_1$ and $m_2$ symmetrically as long as the criteria \eqref{eq:h1} and \eqref{eq:h2} are satisfied after optimization.

\textit{Domain of satisfaction ---} The optimization of loss function \eqref{eq:loss} requires providing it with the values of $x$. According to theorem~\ref{theo:Lyapanov}, the criteria \eqref{eq:lyapanov_psd} and \eqref{eq:lyapanov_nsd} must be satisfied in all points $x$ of a domain $D$ about the origin $\bar{x}=0$. There are two issues about this domain. The first one is that we have no information about the size and shape of $D$. The second problem is that, D is a continuous space with infinite number of elements but the loss function \eqref{eq:loss} can only be evaluated on a discrete set of points. To mitigate these problems, we introduce two parameters $\{\delta, r\}$. Assume we restrict domain $D$ to a ball about the origin and $r$ is the radius of that ball. This assumption is valid because whatever continuous compact set $D$ we consider about $x=0$, we can always construct a ball centering at $x=0$ within $D$ and choose this ball as a new domain $D_r$ about the origin. Of course this choice for $D$ is conservative but assists us by reducing the number of required parameters for characterizing $D$. Once the domain $D_r$ is determined, we need to sample from it to optimize loss function \eqref{eq:loss}. This can be done in stochastic (randomly distributed samples) or deterministic (regularly distributed samples) way. We found out that the stochastic method is much more effective than deterministic one; so here we only explain the stochastic method. 
In this approach, we take random samples from domain $D_r$. Parameter $\delta$ determines how finely we sample from $D$. For instance if the samples are generated uniformly from $D_r$ (see appendix~\ref{sec:sample_sphere}) by resolution $\delta$, the number of samples $N$ is roughly determined by $N=(2r/\delta)^d$(see appendix~\ref{sec:sample_sphere}). In general, assume generated samples come from distribution $p_{D_r}(x;\delta)$ where the subscript $D_r$ shows the domain of interest and the resolution of sampling is controlled by parameter $\delta$. The following loss function is then minimized for $\theta$:

\begin{equation}
\theta^*=\argmin_\theta \expectation_{x\sim p_{D_r}(x;\delta)} \mathcal{L}(x;\theta) 
\label{eq:stochastic_loss}
\end{equation}


Minimizing the loss function \eqref{eq:stochastic_loss} ensures the conditions of theorem~\ref{theo:Lyapanov} when $\mathcal{L}\to 0^+$ and $\delta \to 0$.

If the loss function decreases to 0 over iterations of optimization, $V(x;\theta^*)$ is proposed as a Lyapunov function with desired conditions of theorem~\ref{theo:Lyapanov} and the system $\dot{x}=f(x)$ is asymptotically stable. If the loss function does not converge to 0, we cannot say anything about the stability of the system. However, there is subtle point that is worth mentioning here. Since we are using universal function approximators, we can loosely argue that SGD is searching in the space of \emph{all possible functions}. Therefore, its failure in finding the Lyapunov function suggests \emph{instability} of the system.

\subsection{Experiments}
Here we test \delf on a couple of dynamical systems to show its efficacy in determining the local stability. The description of each experiment comes in the caption of the corresponding tables (~\ref{tbl:2states}, ~\ref{tbl:3states} and ~\ref{tbl:unstable}). See appendix.~\ref{sec:network_architecture} for architectural and optimization details.

\begin{table}[t!]
     \begin{center}
     \begin{tabular}{ c c c c c c c }
     \toprule
\multicolumn{3}{c}{
\text{Stable:}\(\displaystyle
\left\{
	\begin{array}{ll}
	\dot{x}_1=x_1-x_1^3+x_2\\
	\dot{x}_2=3x_1-x_2\\
	\end{array}
\right.
\)
}
&
\multicolumn{3}{c}{
\text{Stable:}\(\displaystyle
\left\{
	\begin{array}{ll}
	\dot{x}_1=x_2\\
	\dot{x}_2=-\frac{g}{l}\sin x_1-\frac{k}{m}x_2 \text{\tiny ($g=l=k=m=1$)}\\
	\end{array}
\right.
\) 
}
\\
\cmidrule(r){1-3}\cmidrule(lr){4-6}
\vspace{-0.5cm}
\\
\raisebox{-\totalheight}{\includegraphics[width=0.15\textwidth]{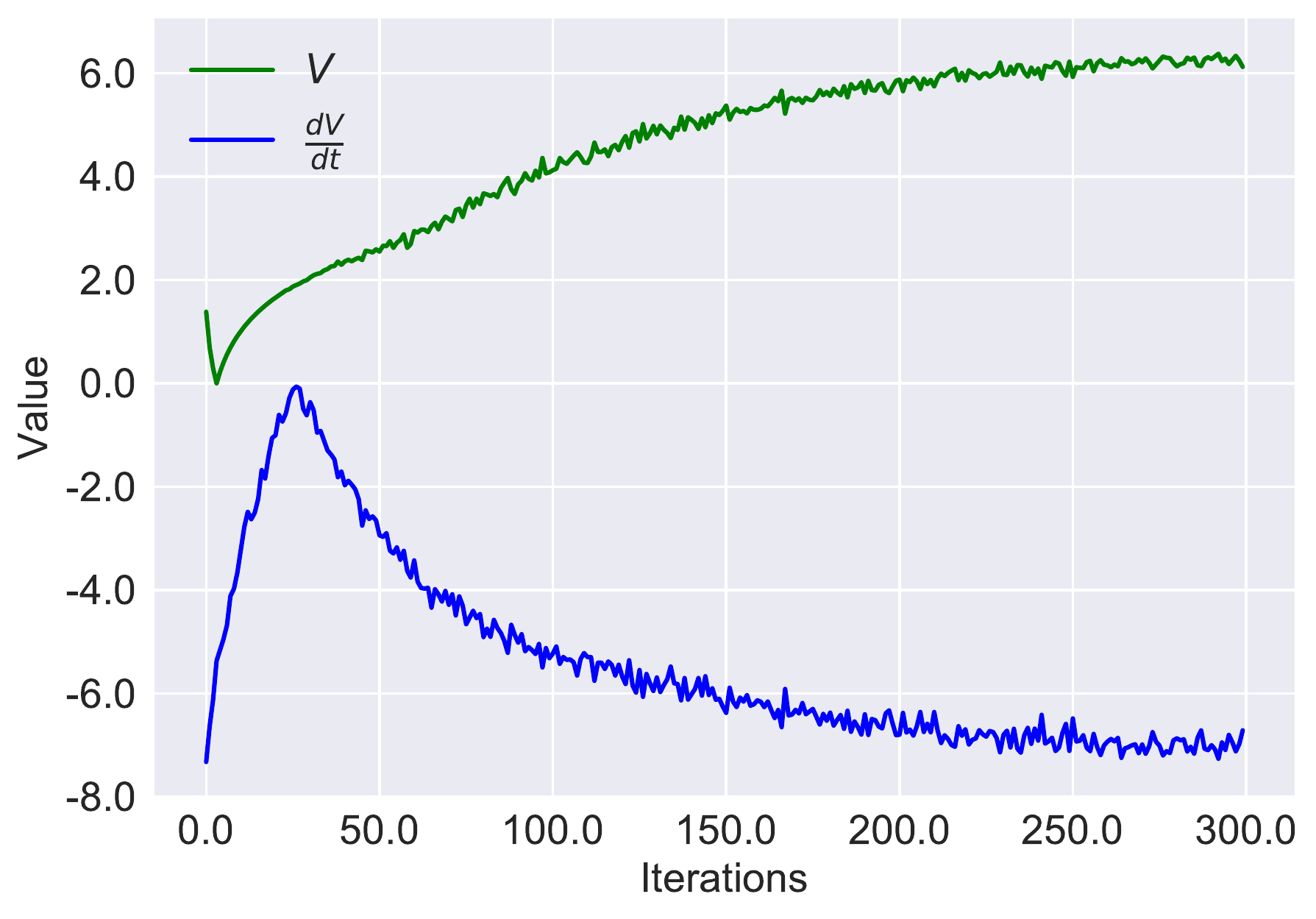}}&
\raisebox{-\totalheight}{\includegraphics[width=0.15\textwidth]{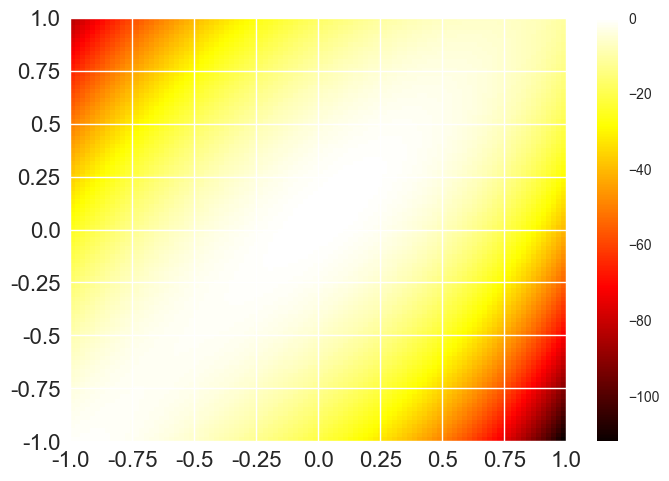}}&
\raisebox{-\totalheight}{\includegraphics[width=0.15\textwidth]{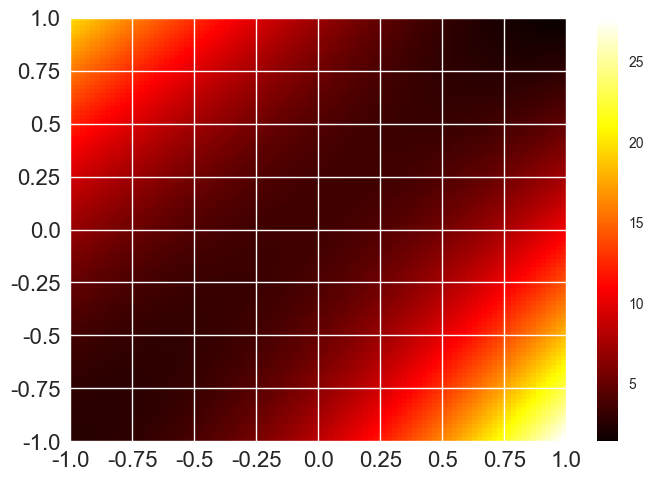}}&
\raisebox{-\totalheight}{\includegraphics[width=0.15\textwidth]{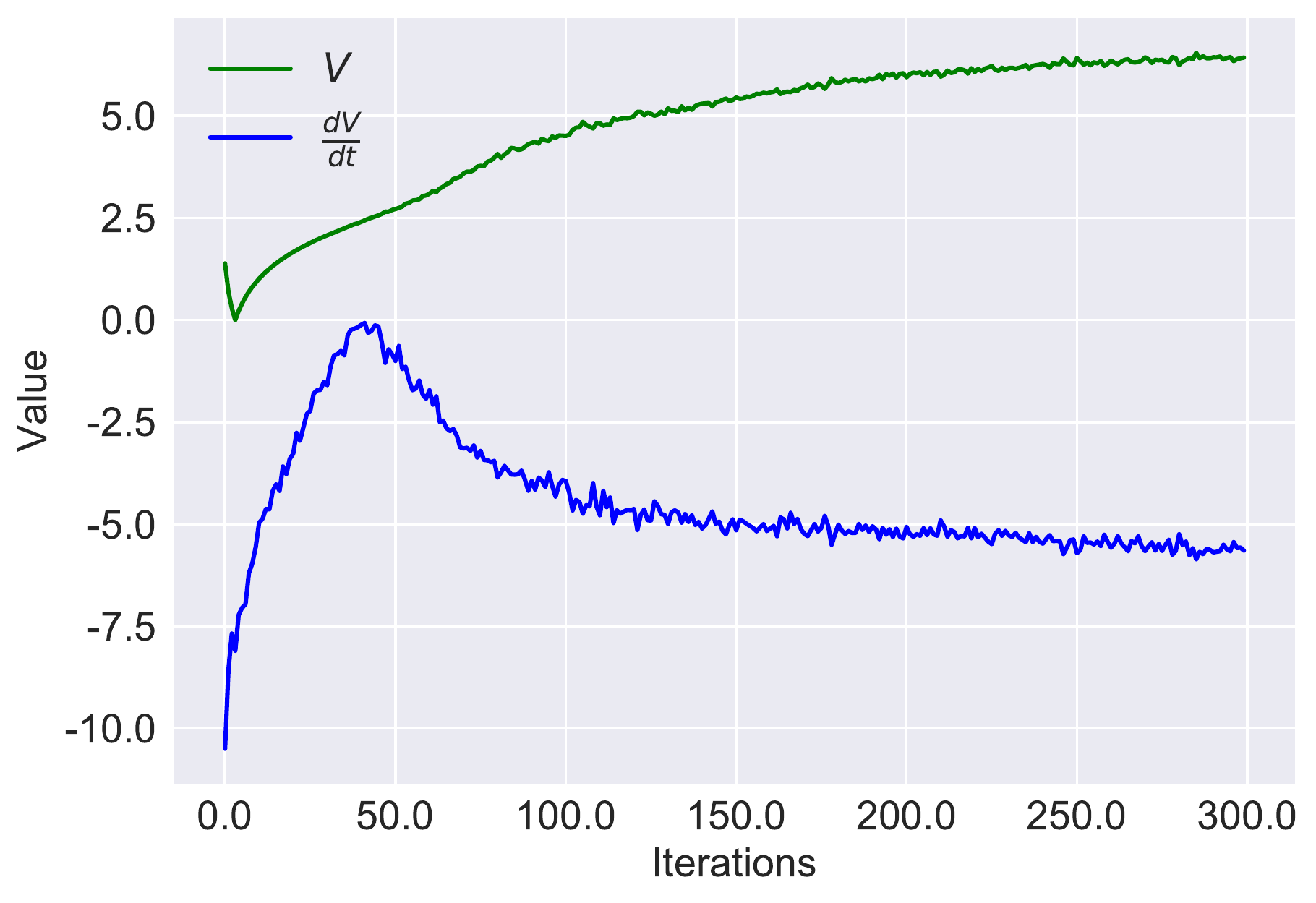}}&
\raisebox{-\totalheight}{\includegraphics[width=0.15\textwidth]{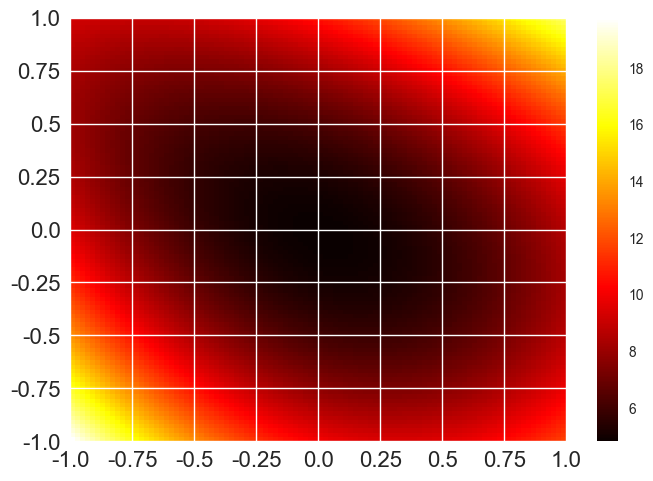}}&
\raisebox{-\totalheight}{\includegraphics[width=0.15\textwidth]{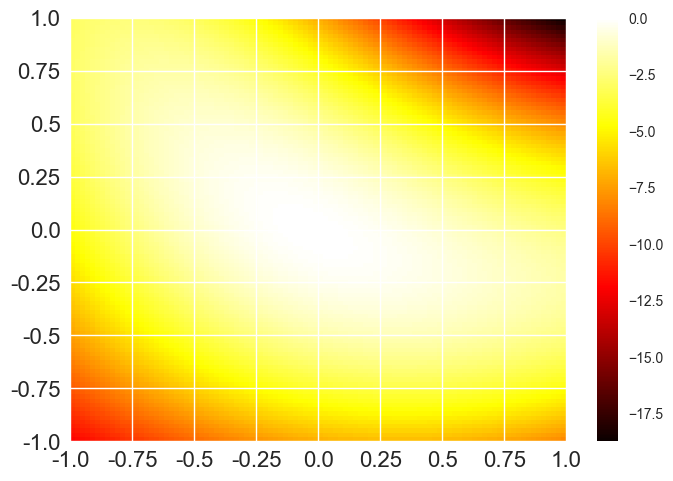}}&
\\
\bottomrule
\end{tabular}
\caption{\small Stability analysis for systems with $2$-dimensional state space. The figures under each system: Left: $\bar{V}=\expectation_{x\sim p_{D_r}(x;\delta)}V(x)$ and $\bar{\dot{V}}=\expectation_{x\sim p_{D_r}(x;\delta)}\dot{V}(x)$, Middle: $V(x)$ around the equilibrium, Right: $\dot{V}(x)$ around the equilibrium}
\label{tbl:2states}
\end{center}
\end{table}

\begin{table}[t!]
     \begin{center}
     \begin{tabular}{ c  p{5cm}  p{5cm}  }
     \toprule

\text{Stable:}\(\displaystyle
\left\{
	\begin{array}{ll}
	\dot{x}_1=-2x_1+x_1^3\\
	\dot{x}_2=-x_2+x_1^2\\
	\dot{x}_3=-x_3
	\end{array}
\right.
\) &

\text{Stable:}\(\displaystyle
\left\{
	\begin{array}{ll}
	\dot{x}_1=-x_1\\
	\dot{x}_2=-x_1-x_3-x_1x_3\\
	\dot{x}_3=(x_1+1)x_2\\
	\end{array}
\right.
\)&
\text{Stable:}\(\displaystyle
\left\{
	\begin{array}{ll}
	\dot{x}_1=-x_2x_3+1\\
	\dot{x}_2=x_1x_3-x_2\\
	\dot{x}_3=x_3^2(1-x_3)
	\end{array}
\right.
\)\\

\cmidrule(r){1-1}\cmidrule(lr){2-2}\cmidrule(l){3-3}
\vspace{-0.5cm}
       \\

\raisebox{-\totalheight}{\includegraphics[width=0.3\textwidth]{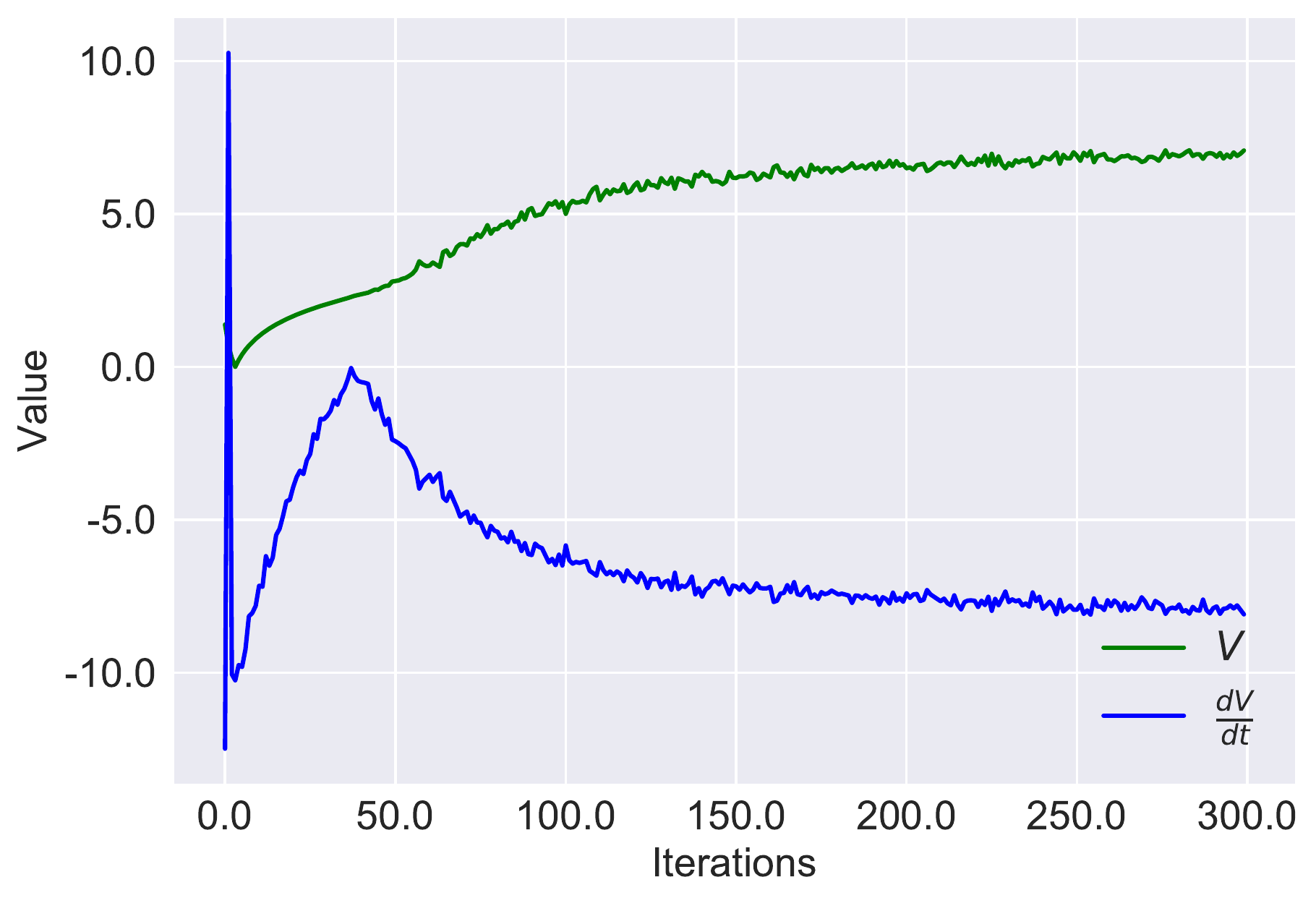}}&
\raisebox{-\totalheight}{\includegraphics[width=0.3\textwidth]{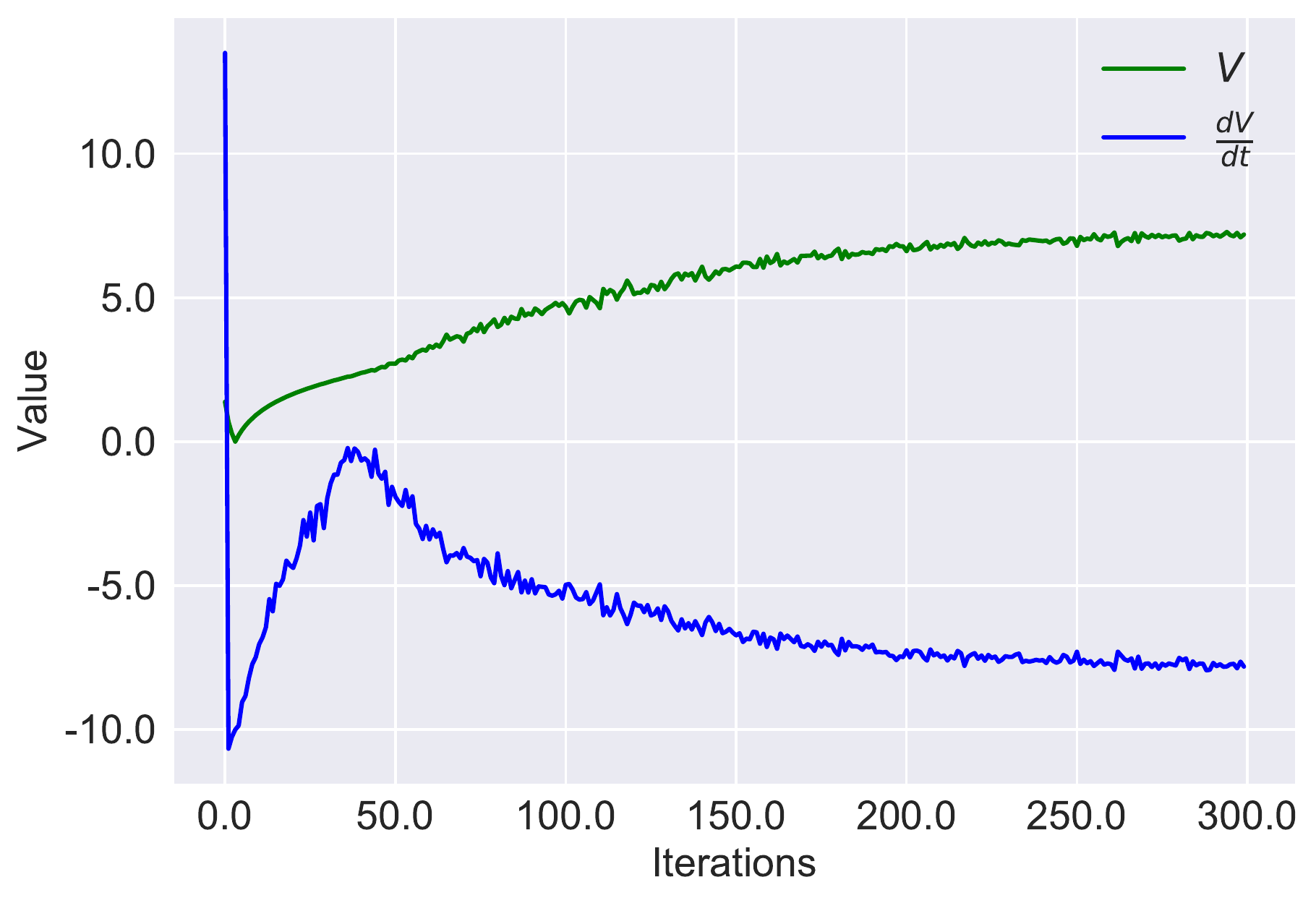}}&
\raisebox{-\totalheight}{\includegraphics[width=0.3\textwidth]{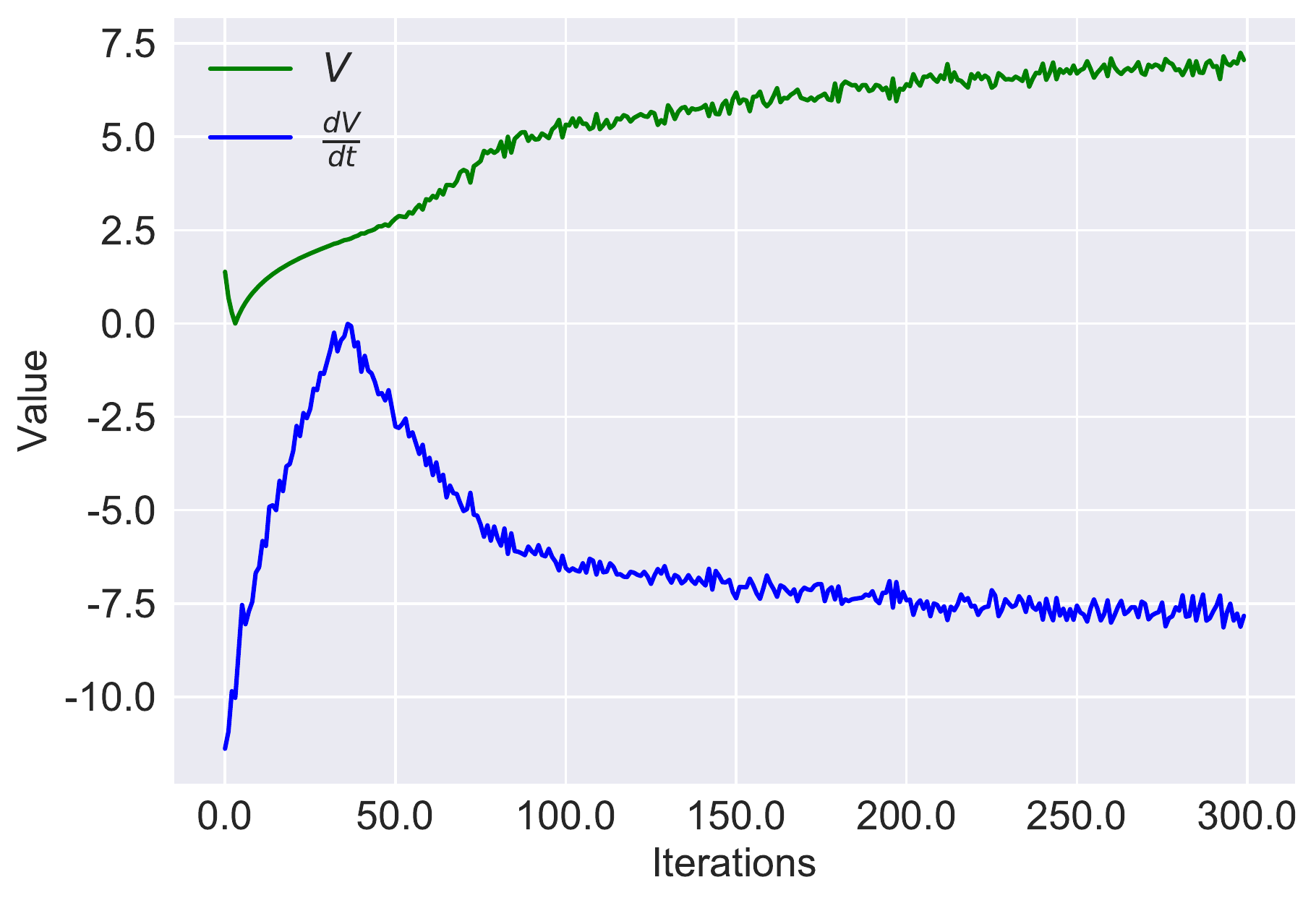}}
      \\ \bottomrule
      \end{tabular}
      \caption{\small Stability analysis for systems with $3$-dimensional state space. The trajectories of the values of $V$ and $\bar{\dot{V}}$ deviates such that $\bar{V}$ becomes persistently positive and $\bar{\dot{V}}$ becomes persistently negative. This pattern of values implies the \textbf{\underline{stability}} of the system.}
      \label{tbl:3states}
      \end{center}
      \end{table}

\begin{table}[t!]
     \begin{center}
     \begin{tabular}{ c  p{5cm}  p{5cm}  }
     \toprule

\(\displaystyle
\text{Unstable:}\left\{
	\begin{array}{ll}
	\dot{x}_1=-x_1 + x_2^2\\
	\dot{x}_2=2x_2-x_1^3\\
	\end{array}
\right.
\) &

\(\displaystyle
\text{Unstable:}\left\{
	\begin{array}{ll}
	\dot{x}_1= x_1 - x_2\\
	\dot{x}_2= -x_1^2 + x_2\\
	\end{array}
\right.
\)&
\(\displaystyle
\text{Unstable:}\left\{
	\begin{array}{ll}
	\dot{x}_1= 3x_1 - x_2\\
	\dot{x}_2= -x_1^3 + 4x^2\\
	\dot{x}_3= x_3
	\end{array}
\right.
\)\\
\cmidrule(r){1-1}\cmidrule(lr){2-2}\cmidrule(l){3-3}
\vspace{-0.5cm}
       \\
\raisebox{-\totalheight}{\includegraphics[width=0.3\textwidth]{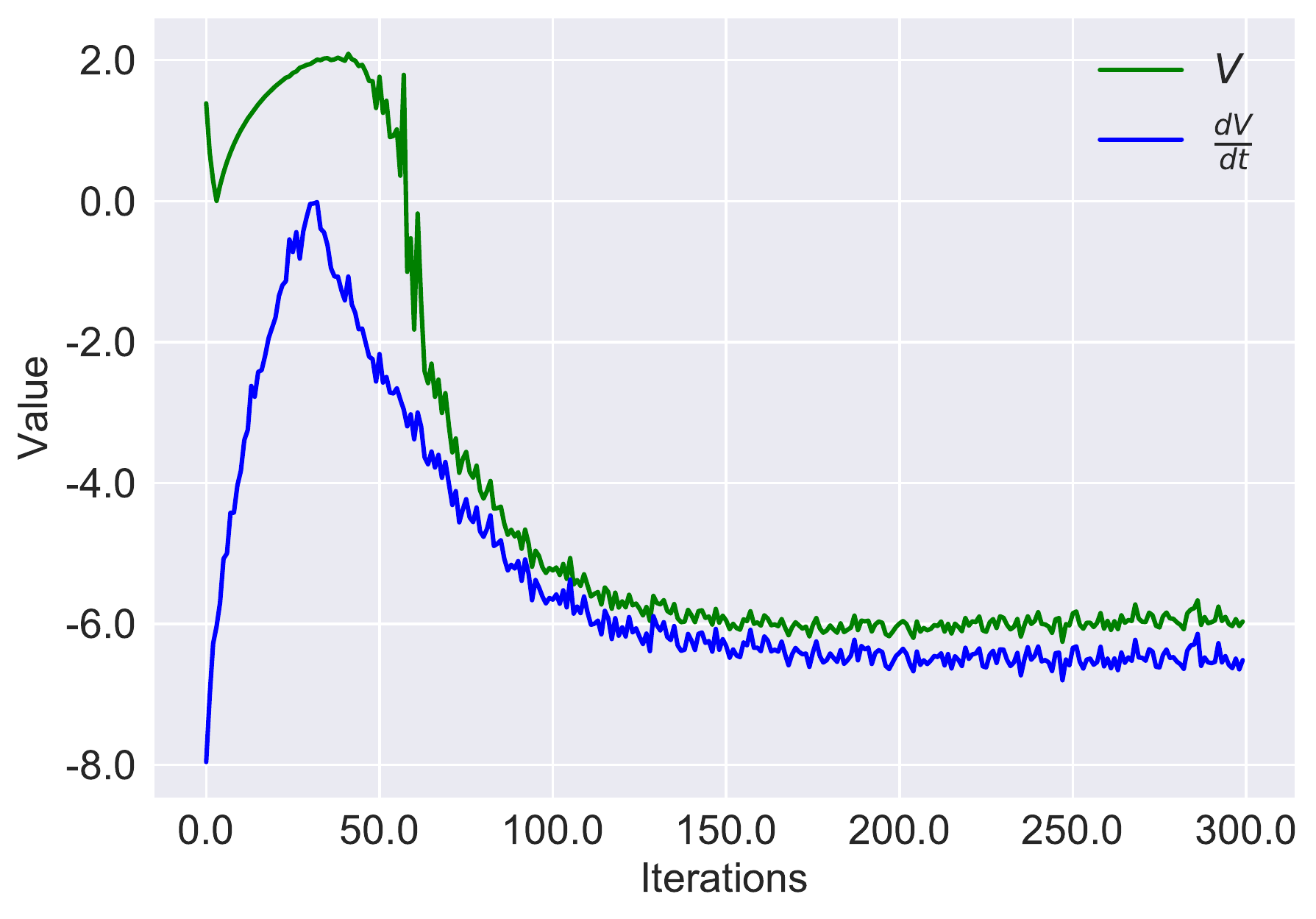}}&
\raisebox{-\totalheight}{\includegraphics[width=0.3\textwidth]{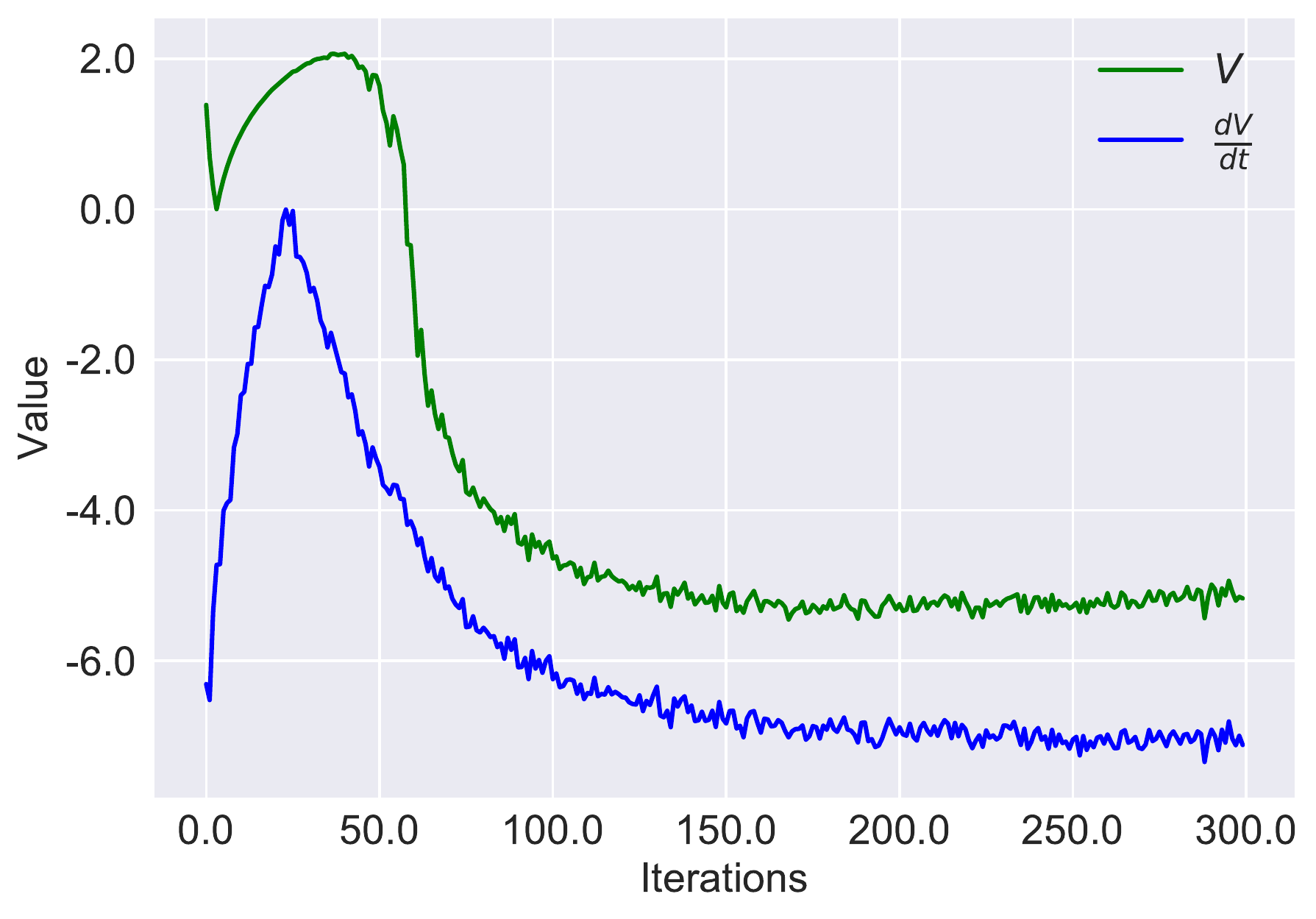}}&
\raisebox{-\totalheight}{\includegraphics[width=0.3\textwidth]{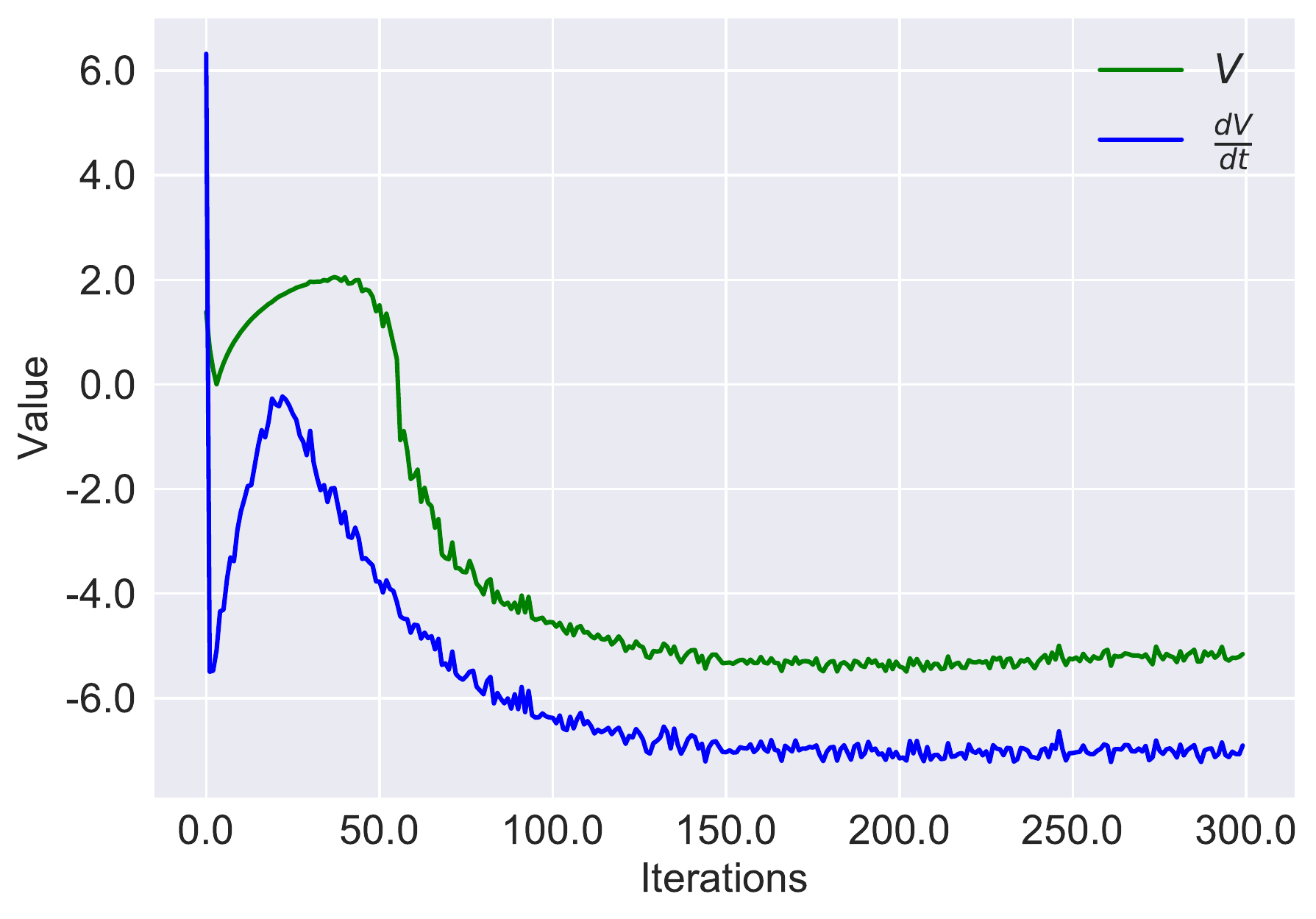}}
      \\ \bottomrule
      \end{tabular}
      \caption{\small The values of $\bar{V}$ and $\bar{\dot{V}}$ during the course of optimization. This pattern is different from the trajectory of the same functions for stable systems (Tables~\ref{tbl:2states} and~\ref{tbl:3states}). In general, we cannot conclude anything about the stability of these systems. But given the universal function approximator assumption for MLPs with sufficient capacity, these trajectory implies the \textbf{\underline{instability}} of the systems since NO Lyapunov function existed in the space of all possible functions.}
      \label{tbl:unstable}
      \end{center}
      \end{table}

\clearpage
\section{Discussion}
In this paper, a simple generic method (\delf) is proposed to investigate the stability of dynamical systems by searching for a Lyapunov function. Because of the stochastic nature of the method, it only provides theoretical guarantee in the limit when the number of samples from the domain $D_r$ goes to infinity (i.e. $N\to\infty$). However, it is still helpful to obtain a quick insight into the stability of dynamical systems whose Lyapunov functions are extremely difficult for human to construct. We are currently working on theoretical guarantee that supports the stability analysis for a Lyapunov function and is learned by a multilayer perceptron.


\bibliographystyle{icml2019}      
\bibliography{main}

\newpage

\section*{Appendices}
\appendix

\section{Drawing samples from within a hypersphere}
\label{sec:sample_sphere}
\begin{figure}[t!]
	\centering
\subfigure[Non-Uniform sampling]{
\includegraphics[width=0.47\linewidth]{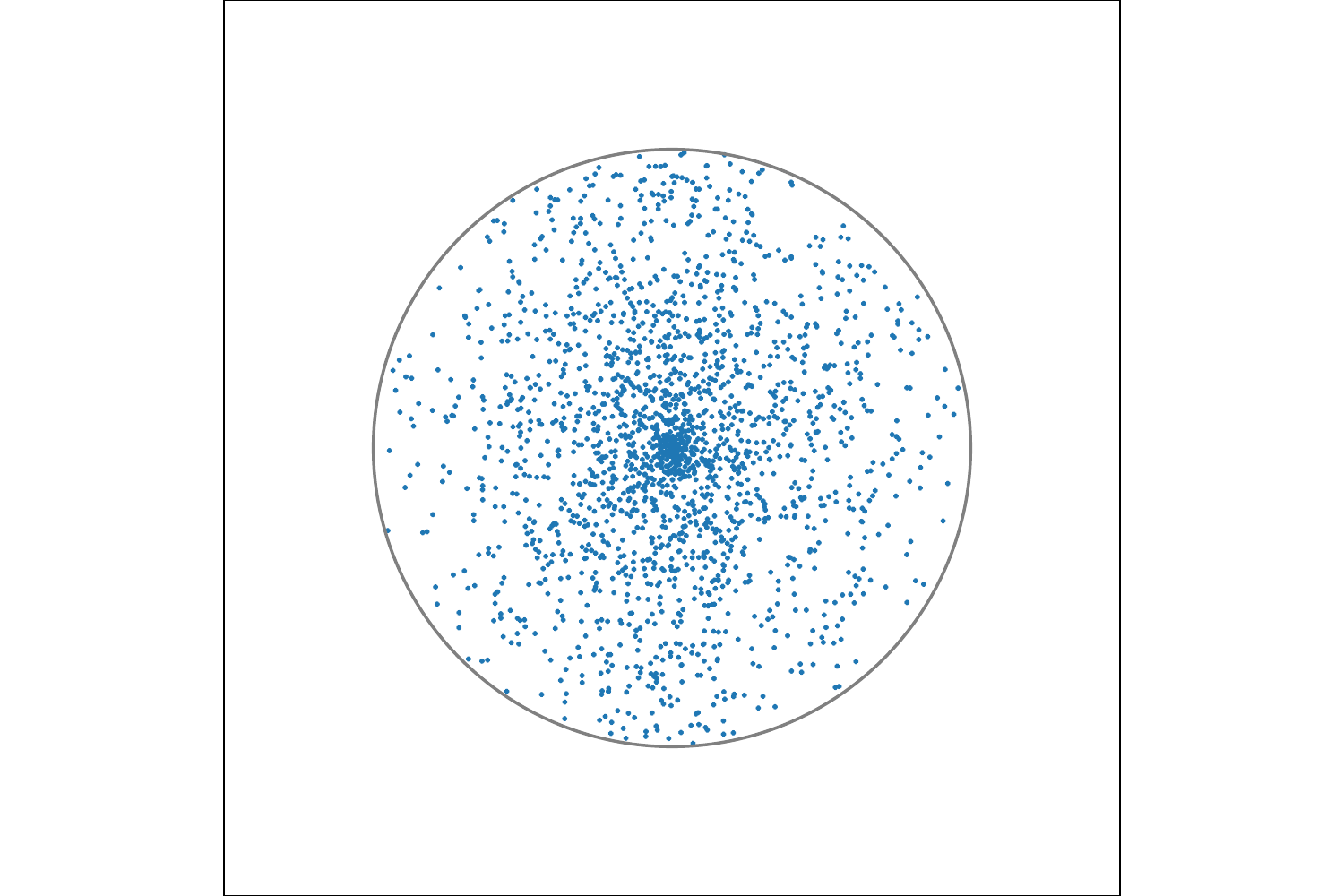}}
\subfigure[Uniform Sampling]{
\includegraphics[width=0.47\linewidth]{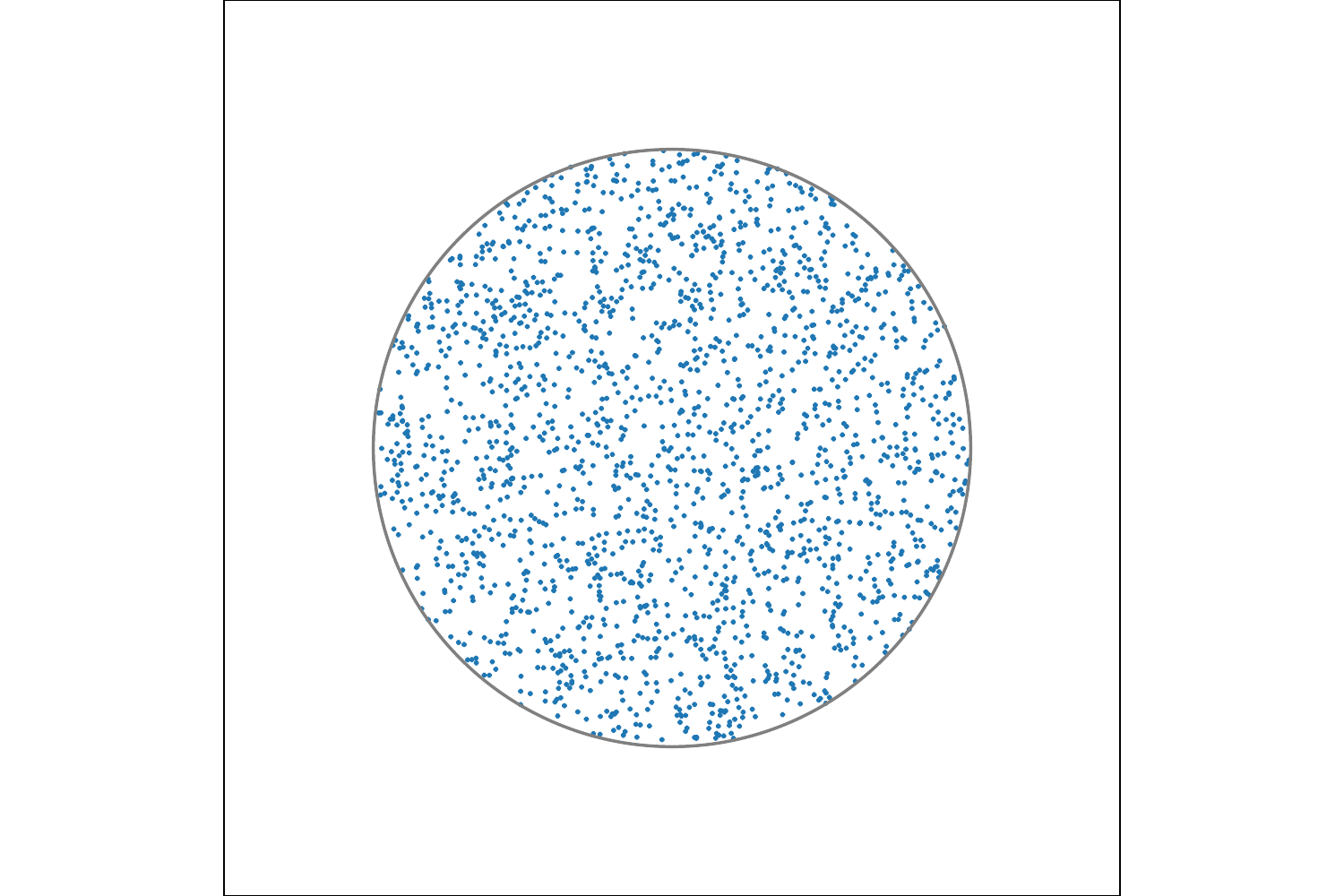}}
\caption{(a) Sampling from a hypersphere with higher resolution around the center. (b) Uniform sampling from a hypersphere }
\label{fig:sampling_sphere}
\end{figure}

We tested two methods for drawing samples from within a hypersphere $D_r=\{x\in \mathbb{R}^n:\lVert x\rVert\leq r\}$:

{\it Center-concentrated sampling--- } In this method, we parameterize the space within the sphere by a polar coordinate in $n$ dimensions. The parameterization consists of a radius $r$ and $n-1$ angles. The relationship between polar coordinates and euclidean coordinates in $n$ dimensions is simply the extension of $3$-dimensional $\{r, \theta, \phi\}$ as follows:
\begin{align}
\label{eq:polar_to_euclidean}
    x_1&=r\cos(\phi_1)\\\nonumber
    x_2&=r\sin(\phi_1)\cos(\phi_2)\\\nonumber
    x_3&=r\sin(\phi_1)\sin(\phi_2)\cos(\phi_3)\\\nonumber
    \vdots\\\nonumber
    x_{n-1}&=r\sin(\phi_1)\ldots\sin(\phi_{n-2})\cos(\phi_{n-1})\\\nonumber
    x_n&=r\sin(\phi_1)\ldots\sin(\phi_{n-2})\sin(\phi_{n-1})\nonumber
\end{align}
We draw random samples uniformly from each of the polar coordinates in their feasible range and map the samples back to the euclidean coordinates via ~\eqref{eq:polar_to_euclidean}. It might look counter-intuitive that the generated samples are distributed more densely around the center even thought each polar coordinate was sampled uniformly (see Fig.~\ref{fig:sampling_sphere}(a)). This is caused by the above nonlinear transformations. Our motivation for using this type of sampling was the definition~\eqref{def:stablility} of stability. To prove the stability in the sense of~\eqref{eq:def_stable}, it is enough to show the existence of a hypersphere with radius $\epsilon$. Therefore, it makes sense to focus our attention on areas closer to the equilibrium. This means that area near the equilibrium is sampled more frequently and forms the major portion of the loss function ~\eqref{eq:loss}.

{\it Uniform sampling--- } In this method, we draw samples which are uniformly distributed all over the hypersphere (see Fig.~\ref{fig:sampling_sphere}(b)). To this end, the following steps must be taken:
\begin{itemize}
    \item Generate $N$ samples from an $n$-dimensional Gaussian distribution.\\
    \begin{equation}
        X_i\sim\mathcal{N}(\mathbf{0}, \mathbf{1})
    \end{equation}
    \item Compute the sum of squared of the normal variables to achieve $Q\sim\chi^2$-squared distribution with $N$ degrees of freedom\\
    \begin{equation}
        Q = \sum_{i=1}^{N}X_i^2
    \end{equation}
    \item Apply the cumulative distribution function of $\chi^2$-squared distributed random variables on the samples of the previous step. This function is called \emph{incomplete gamma} which is defined as
    \begin{equation}
        \Gamma(s, x)=\int_x^\infty t^{s-1}e^{-t}\text{d}t
    \end{equation}
    where $s = N/2$.
    \item  Transform these samples affinely by the center and radius of the hypersphere so that the resultant samples uniformly cover the hypersphere.
\end{itemize}

In our tests, both methods gave comparable results but the uniform sampling method gave faster rate of convergence for functions $V$ and $\dot{V}$. Therefore, we stuck with the uniform sampling in the experiments. It is predictable that the \emph{center-concentrated} sampling can be useful for dynamical systems with more complicated dynamics around the equilibrium.

\section{Sampling resolution}
We use the hyper-parameter $\delta$ to control the number of samples which are needed to be uniformly drawn from a hypersphere $D_r$ with radius $r$ in $d$-dimensional space. According to the definition of $\delta$, the average distance between two samples in $d$-dimensional space must be at least $\delta$ meaning that each sample is located at the center of a tiny hypersphere $D_{\delta/2}$ with radius $\delta/2$. When $N$ samples are drawn uniformly from within $D_r$, we can assume $D_r$ is filled with $N$ tiny hyperspheres $D_{\delta/2}$. This implies the equality of two volumes in $d$-dimensional space:
\begin{equation}
    \frac{\pi^\frac{d}{2}}{\Gamma(\frac{d}{2}+1)}r^d \approx N\frac{\pi^\frac{d}{2}}{\Gamma(\frac{d}{2}+1)}(\frac{\delta}{2})^d\implies N\approx 2^d(\frac{N}{\delta})^d
\end{equation}
where $\Gamma(.)$ is the gamma function and the left-hand side is the volume of a $d$-dimensional hypersphere with radius $r$.

\begin{figure}[t!]
	\centering
\subfigure[Physical System]{
\includegraphics[width=0.37\linewidth]{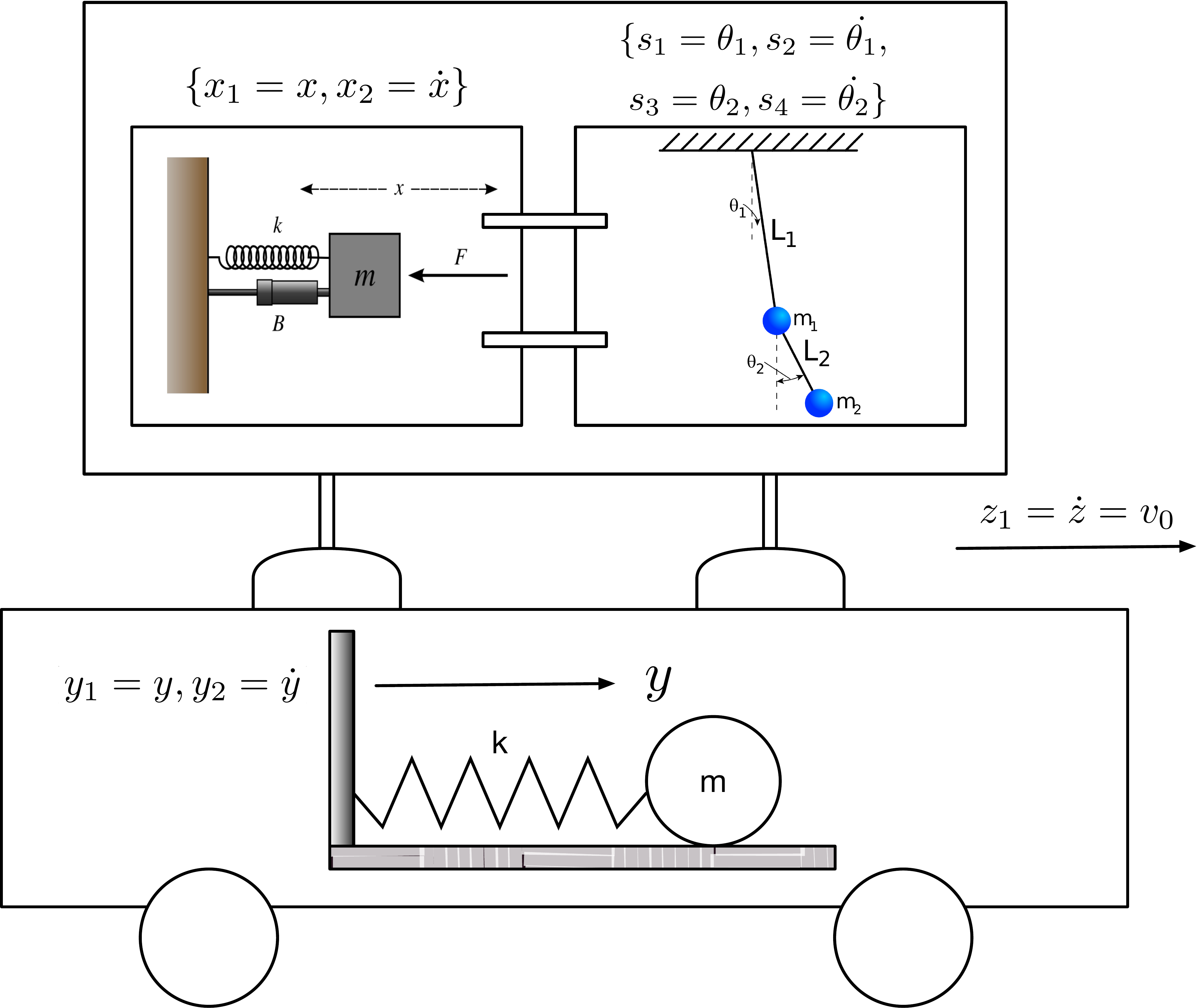}}
\hspace{1cm}
\subfigure[Energy Hierarchy]{
\includegraphics[width=0.37\linewidth]{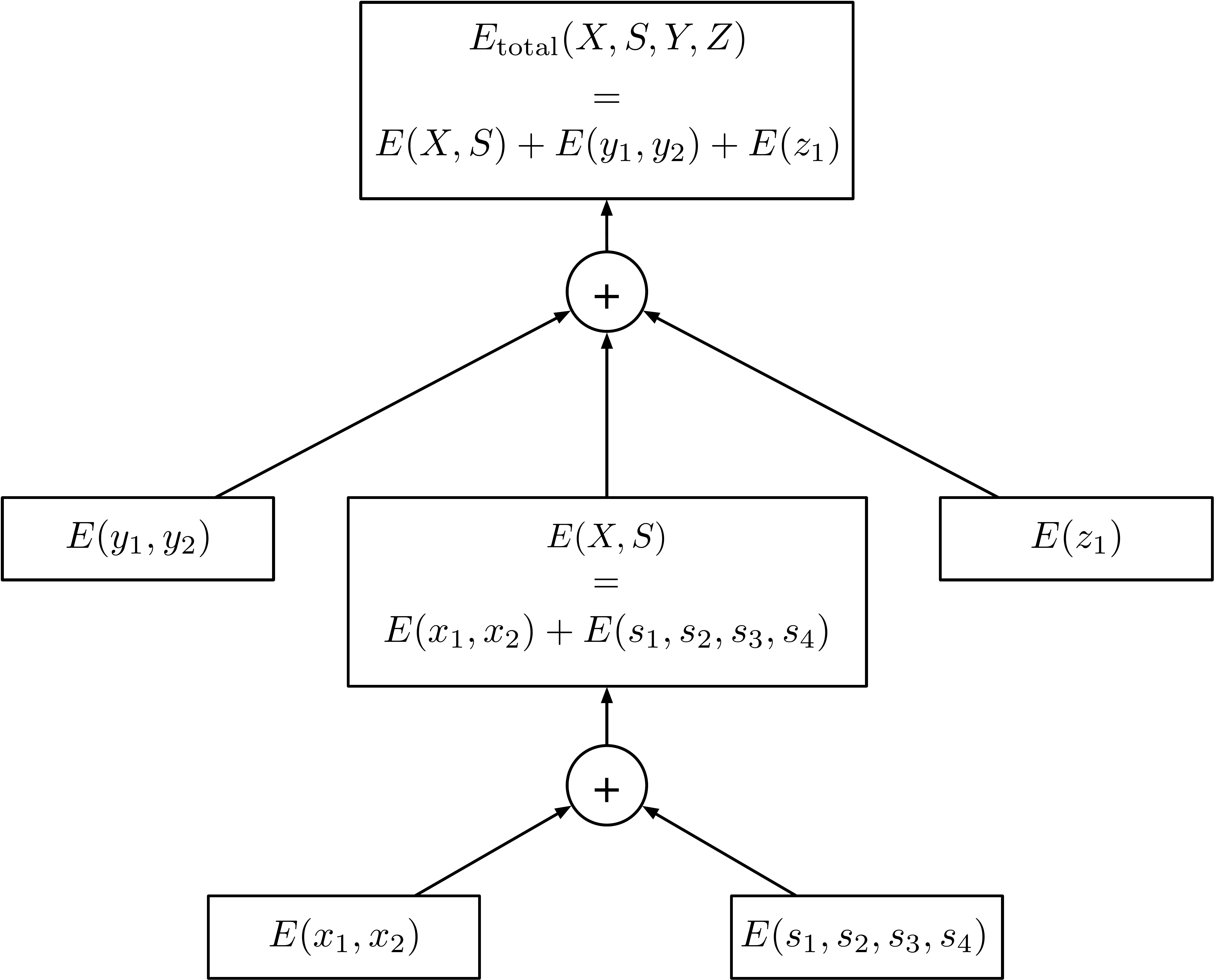}}
\caption{Left: A hypothetical physical system with several internal dynamics. The total energy of the system can be seen as the sum of the energy of each of these internal dynamics. Right: The hierarchy of the energies of involved internal dynamics that constitute the total energy of the system.}
\label{fig:energy_schematics}
\end{figure}
\section{Details on the experiments}
\label{sec:network_architecture}
Neural networks with sufficiently wide hidden layer can approximate any function assuming some weak conditions (smoothness, etc)~\citep{csaji2001approximation}. However, it is always helpful to incorporate inductive bias as a guide to ease solving the problem. Lyapunov function has close relationship with energy in physical systems. On the other hand, it is reasonable to assume a hierarchical structure for energy in physical systems. This come from the additive nature of energy and also possibility to find separate clusters of tightly interacting variables belonging to almost separate subsystems~\citep{chen1998linear}. This justifies the use of a multilayer neural network for modeling the Lyapunov function. The number of layers can roughly express our belief about the number of steps in the energy hierarchy. To make it more clear, the hypothetical system of Fig.~\ref{fig:energy_schematics}(a) shows the idea of hierarchical energy in physical systems. In this case, a three layer neural network will have the potential to recover this hierarchy. Notice that this inductive bias is neither necessary nor sufficient for finding a good Lyapunov function. However, it can be a good starting point. Another inductive bias is the prevalence of polynomials in Lyapunov functions. Many of the physical and abstract dynamical systems have at least one polynomial Lyapunov function~\citep{khalil1996noninear}. This observation suggests using polynomials as the activation function for the MLP that approximate \delf. For the experiments, we used a 3 layer neural network with hidden dimension 5 and polynomials of degrees 2 and 3 as activation functions. We also observed that a network with thirds order polynomial followed by a linear layer works well for most cases. This is expected since it is known that many of these dynamical systems have Lyapunov functions which is second or third order polynomial of states. Also, it was observed that three layer smooth MLP with tanh nonlinearity worked well for all cases. This is also expected since this network can approximate any polynomial locally around the equilibrium point. In all experiments, we used SGD as the optimizer with batch size 20 and learning rate 0.005.

\end{document}